\documentclass[11pt]{article}
\usepackage{amsfonts}
\usepackage{latexsym, amssymb, amsmath, amscd, amsfonts, epsfig, graphicx, colordvi}
\usepackage{ifpdf}
\usepackage{blkarray}
\usepackage{graphicx}
\usepackage{xcolor}
\usepackage[dvipsnames]{pstricks}
\newtheorem{thm}{Theorem}[section]

\newtheorem{defi}[thm]{Definition}

\makeatletter
\def\ExtendSymbol#1#2#3#4#5{\ext@arrow 0099{\arrowfill@#1#2#3}{#4}{#5}}
\makeatother

\title{Combinatorial Proofs of Two Overpartition Theorems Connected by a Universal Mock Theta Function}
\author{
Doris D. M. Sang\raisebox{5pt}{\scriptsize 1} and Diane Y. H.
Shi\raisebox{5pt}{\scriptsize 2}}
\date{
\vspace{15pt}\raisebox{5pt}{\scriptsize 1\,}School of Mathematics
and Quantitative Economics\\Dongbei
University of Finance and Economics, Liaoning 116025, P.R. China\\sdm@cfc.nankai.edu.cn\\
\vspace{15pt}\raisebox{5pt}{\scriptsize 2\,}School of
Mathematics\\ Tianjin University, Tianjin 300072, P.R.
China\\shiyahui@tju.edu.cn}

\begin{document}

\maketitle

\vspace{0.3cm} \noindent{\bf Abstract.} In 2015,  Bringmann, Lovejoy and Mahlburg  considered certain kinds overpartitions, which can been seen as the overpartition analogue of Schur's partition. The motivation of their work is that  the difference   between the generating function of Schur's classical partitions and the generating functions of the partitions in which the smallest part is excluded. The difference between the two generating functions of partitions is a specialization of the ¡°universal¡± mock theta function $g_3$ which introduced by Hickerson. To give an analogue of this, by using another universal mock theta function $g_2$ instead of $g_3$, Bringmann Lovejoy and Mahlburg introduced  two kinds of overpartitions, which satisfy certain congruence conditions and difference conditions with the smallest parts different. They prove these theorems by using the $q$-differential equations. In this paper, we will give the generating functions of these two kinds of overpartitions by combinatorial technique.



\section{Introduction}

In 2015,  Bringmann, Lovejoy and Mahlburg \cite{bri15}  considered some overpartition  theorems,  motivated by the relation  between the generating function of Schur's classical partitions and the generating functions of the partitions in which the smallest parts are excluded. The difference between the two generating functions of partitions is a specialization of the ¡°universal¡± mock theta function $g_3(x;q)$ up to the odd order. To give an analogue of this, by using $g_2(x;q)$ instead of $g_3(x;q)$, Bringmann, Lovejoy and Mahlburg \cite{bri15} introduced  two kinds of overpartitions, which satisfy certain congruence conditions and difference conditions with different smallest parts. They proved these theorem by using the $q$-differential equations. At the end of \cite{bri15}, the authors asked for the combinatorial interpretation of these theorems. In this paper, we will give the generating functions of these two kinds of overpartitions by combinatorial method.

Let us give an overview of  some definitions. A partition $\lambda$
of a positive integer $n$ is a non-increasing sequence of positive
integers $\lambda_1\geq \cdots\geq \lambda_s>0$ such that
$n=\lambda_1+\cdots+\lambda_s$. The partition of zero is the
partition with no parts. An overpartition $\lambda$ of a positive
integer $n$ is also a non-increasing sequence of positive integers
$\lambda_1\geq \cdots\geq \lambda_s>0$ such that
$n=\lambda_1+\cdots+\lambda_s$ and the first occurrence  of  each
integer may be overlined. For example,
$(\overline{7},7,6,\overline{5},2,\overline{1})$ is an overpartition
of $28$.

We shall adopt  the common notation  as used in Andrews \cite{and76}. Let
 \[(a)_\infty=(a;q)_\infty=\prod_{i=0}^{\infty}(1-aq^i),\]
 and \[(a)_n=(a;q)_n=\frac{(a)_\infty}{(aq^n)_\infty}.\]
We also write
 \[(a_1,\ldots,a_k;q)_\infty=(a_1;q)_\infty\cdots(a_k;q)_\infty.\]

We first recall the following classical Schur's partition theorem \cite{sch26}. Throughout the paper we assume that $d\geq 3$ and
$1 \leq r<d/2$.
\begin{thm}[Schur]\label{Schur}For all $n\geq 0$, let $B_{d,r}(n)$ denote the number of partitions of n into parts congruent
to $r$, $d-r$, or $d \pmod{d}$ such that $\lambda_{i+1}-\lambda_{i}\geq d$ with strict inequality if $d|\lambda_{i+1}$. Let $E_{d,r}(n)$
denote the number of partitions of $n$ into distinct parts that are congruent to $\pm r \pmod{d}$. Then
\[E_{d,r}(n)=B_{d,r}(n).\]
\end{thm}
Bringmann and Mahlburg \cite{bri14} denoted the generating function for $B_{d,r}(n)$ by
\[\mathfrak{B}_{d,r}(q) :=\sum_{n\geq 0}B_{d,r}(n)q^n.\]
 Then
 \begin{equation}\label{schur}\mathfrak{B}_{d,r}(q)=\sum_{n\geq 0}E_{d,r}(n)q^n=(-q^r,-q^{d-r};q^d)_\infty,\end{equation}
which implies that $\mathfrak{B}_{d,r}(q)$ is a modular function (up to a rational power of $q$).

By giving an additional restriction on the smallest parts, Bringmann and Mahlburg \cite{bri14} defined a new kind of partitions, whose generating function  equals the product of $\mathfrak{B}_{d,r}(q)$ and a certain
specialization of a  universal mock theta function $g_3(x;q)$. This result is stated as follows.
\begin{thm}
Let $C_{d,r}(n)$ denote the number of partitions enumerated by $B_{d,r}(n)$ that also satisfy the
additional restriction that the smallest part is larger than $d$. Then the generating function for $C_{d,r}(n)$ is
\begin{equation}\label{cg3}\mathfrak{C}_{d,r}(q) :=\sum_{n\geq 0}C_{d,r}(n)q^n
=\mathfrak{B}_{d,r}(q)g_3(-q^r;q^d)=(-q^r,-q^{d-r};q^d)_\infty \sum_{n\geq 0}\frac{q^{dn(n+1)}}{(-q^r,-q^{d-r};q^d)_{n+1}},\end{equation}
where
\[g_3(x; q) :=\sum_{n\geq 0}\frac{q^{n(n+1)}}{(x,q/x;q)_{n+1}}.\]\end{thm}
Here $g_3(x; q)$ is Hickerson's \cite{hic88,hic88b} universal mock theta function (of odd order).
This implies  that
$\mathfrak{C}_{d,r}(q)$ is not a modular form, but the product of the modular form $\mathfrak{B}_{d,r}(q)$ and a
specialization of the universal mock theta function $g_3$ instead.

The universal mock theta function $g_3(x; q)$ is so-named because Hickerson \cite{hic88,hic88b}, Gordon
and McIntosh \cite{gor} have shown that each of the classical odd order mock theta functions may be
expressed, up to the addition of a modular form, as a specialization of $g_3(x; q)$.

There is a second universal mock theta function,
\[g_2(x; q) :=\sum_{n\geq 0}\frac{(-q;q)_nq^{n(n+1)/2}}{(x,q/x;q)_{n+1}},\]
which corresponds to the classical even order mock theta functions \cite{gor}.

To search for an
analogue of \eqref{cg3} with $g_2(x; q)$ in place of $g_3(x; q)$, Bringmann, Lovejoy and Mahlburg \cite{bri15} introduced  two kinds of overpartitions. For the difference conditions between the successive parts are much more complicated, before stated the definitions
 of the two kinds overpartitions, Bringmann, Lovejoy and Mahlburg \cite{bri15} defined a matrix which is indexed by the congruences when modulo $d$.
\begin{defi}Define the $4\times 4$ matrix $\overline{A}_{d,r}$ by
\begin{equation}\label{matrix}
\overline{A}_{d,r}=
\begin{blockarray}{ccccc}
& \overline{r} & \overline{d-r} & \overline{d} & d \\
\begin{block}{c(cccc)}
  \overline{r} & d & 2r & d+r & r           \\
  \overline{d-r} & 2d-2r & d & 2d-r & d-r   \\
  \overline{d} & 2d-r & d+r & 2d & d        \\
  d & d-r & r & d & 0            \\
\end{block}
\end{blockarray}.
 \end{equation}
The rows and columns are indexed by $\overline{r}$, $\overline{d-r}$, $\overline{d}$, and $d$, so that, for example, $\overline{A}_{d,r}(\overline{d}, \overline{d-r})=d+r$.
\end{defi}

By using this matrix, Bringmann Lovejoy and Mahlburg \cite{bri15} introduced the overpartition analogue of Schur's partition function.
\begin{thm}\label{thmB}
For $n\geq 2$, let $\overline{B}_{d,r}(n)$ denote the number of such overpartitions $\lambda$ of
$n$ whose  parts congruent to $r$, $d-r$, or $d \pmod{d}$ and  only multiples of
$d$ may appear non-overlined. In addition, $\lambda$ satisfies the following conditions
\begin{itemize}
\item[(i)] the smallest part is $\overline{r}$, $\overline{d-r}$, $\overline{d}$, or $2d$ modulo $2d$;
\item[(ii)] for $u, v \in \{\overline{r}, \overline{d-r}, \overline{d}, d\}$, if $\lambda_{i+1}\equiv  u \pmod{d}$ and $\lambda_i \equiv v \pmod{d}$, then $\lambda_{i+1}-\lambda_i\geq
\overline{A}_{d,r}(u, v)$;
\item[(iii)] for $u, v \in \{\overline{r}, \overline{d-r}, \overline{d}, d\}$, if $\lambda_{i+1} \equiv u \pmod{d}$ and $\lambda_i\equiv v \pmod{d}$, then $\lambda_{i+1}-\lambda_i\equiv\overline{A}_{d,r}(u, v) \pmod{2d}$. In another word, the actual difference between two parts must be congruent
modulo $2d$ to the smallest allowable difference.
\end{itemize}

Let $\overline{E}_{d,r}(n)$ denote the number of partitions of $n$ into distinct parts congruent to
$\pm r \pmod{d}$ and unrestricted parts divisible by $2d$. Then for all $n\geq  0$, we have

\begin{equation}\label{BE}\overline{B}_{d,r}(n) =\overline{E}_{d,r}(n).\end{equation}
\end{thm}
Denote the generating function for $\overline{B}_{d,r}(n)$ by
\[\overline{\mathfrak{B}}_{d,r}(q):=\sum_{n\geq 0}\overline{B}_{d,r}(n)q^n.\]
Then they  derived  that
\begin{equation}\label{gnrB}\overline{\mathfrak{B}}_{d,r}(q)
=\frac{(-q^r,-q^{d-r};q^d)_\infty}{(q^{2d};q^{2d})_\infty}.\end{equation}

To give an overpartition analogue of Theorem \ref{thmC}, Bringmann, Lovejoy and Mahlburg \cite{bri15} defined another kind of overpartitions enumerated by $\overline{C}_{d,r}(n)$ which have the  different restrictions on the smallest parts. They also gave the generating function of $\overline{C}_{d,r}(n)$, which is the product of the generating function of $\overline{B}_{d,r}(n)$ and  $g_2(-q^r;q^d)$.
\begin{thm}\label{thmC}
Let  $\overline{C}_{d,r}(n)$  be the number of overpartitions of $n$ satisfying conditions (ii) and (iii)
in the definition of $\overline{B}_{d,r}(n)$, with condition (i) modified to be that the smallest part is congruent
to $d$, $\overline{d + r}$, $\overline{2d-r}$ or $\overline{2d} \pmod{2d}$.

Denote the generating
function for $\overline{C}_{d,r}(n)$ by
\[\overline{\mathfrak{C}}_{d,r}(q):=\sum_{n\geq 0}\overline{C}_{d,r}(n)q^n.\]
Then \begin{equation}\label{bc}\overline{\mathfrak{C}}_{d,r}(q) = \overline{\mathfrak{B}}_{d,r}(q)\times g_2(-q^r;q^d)=
\frac{(-q^r,-q^{d-r};q^d)_\infty}{(q^{2d};q^{2d})_\infty}\sum_{n\geq 0}\frac{(-q^d;q^d)_nq^{\frac{dn(n+1)}{2}}}{(-q^r,-q^{d-r};q^d)_{n+1}}.\end{equation}
\end{thm}

Bringmann, Lovejoy and Mahlburg \cite{bri15} proved Theorem \ref{thmB} and \ref{thmC} by  deriving and solving $q$-difference equations
satisfied by the generating functions for the relevant overpartitions.
At the end of \cite{bri15}, they posed a problem that is ``it would be interesting to see a bijective proof of Theorem \ref{thmB}..."

In this paper, we will give an bijective proof of Theorem \ref{thmB}, and derive, by a constructive way, an alternative form of the genrating function of $\overline{C}_{d,r}(n)$ which is equal to the right-hand side of \eqref{bc}. For the difference conditions described by the $4\times 4$ matrix \eqref{matrix} are  complicated,   we employ the $d$-modular Ferrers diagram and $d$-modular partition (overpartition) in the next section, which makes the differences matrix look much easier. In section 3 and 4, we give the combinatorial proof of Theorem \ref{thmB} and Theorem \ref{thmC}.

\section{The $d$-modular Ferrers graph and $d$-modular overpartition}

We shall employ the $d$-modular Ferrers graph to describe the overpartitions  enumerated by $\overline{B}_{d,r}(n)$. By this diagram, we give a new matrix which describes the difference between two successive rows in $d$-modular Ferrers diagram instead of the difference between two parts in $\lambda$. We will rewrite the difference condition matrix $\overline{A}_{d,r}$, which makes the difference conditions   much easier to describe.

\begin{defi}\label{label}For an overpartition $\lambda$ whose overlined parts are $\equiv \pm r,0 \pmod{d}$, non-overlined parts are multiples of $d$, we give the $d$-modular Ferrers graph as follows.
\begin{itemize}
\item[(i)]If the overlined part $\lambda_i=(n-1)d+r\equiv r\pmod{d}$, then we denote it by $(n-1)$ $d$'s and an $\overline{r}$.
\item[(ii)]If the overlined part $\lambda_i=(n-1)d+(d-r)\equiv d-r\pmod{d}$, then we denote it by $(n-1)$ $d$'s and an $\overline{d-r}$.
\item[(iii)]If the overlined part $\lambda_i=nd\equiv 0\pmod{d}$, then we denote it by $(n-1)$ $d$'s and an $\overline{d}$.
\item[(iv)]If the non-overlined part $\lambda_i=(n-1)d\equiv 0\pmod{d}$, then we denote it by $(n-1)$ $d$'s and a $0$.
\end{itemize}
\end{defi}
For example, let $d=7,r=2$, $\lambda=(\overline{26},\ 21,\ 21,\ \overline{16},\ \overline{7})$, then the $7$-modular Ferrers graph of $\lambda$ is as follows:
\begin{center}
 \begin{picture}(100,100)
\put(0,15){$\overline{7}$}

\put(0,30){$7$}\put(15,30){$7$}\put(30,30){$\overline{2}$}

\put(15,45){$7$}
\put(0,45){$7$}
\put(30,45){$7$}\put(45,45){$0$}

\put(15,60){$7$}
\put(0,60){$7$}
\put(30,60){$7$}\put(45,60){$0$}

\put(0,75){$7$}
\put(15,75){$7$}
\put(30,75){$7$}
\put(45,75){$\overline{5}$}
 \end{picture}
\end{center}
According to the $d$-modular Ferrers diagram of $\lambda$, we define a corresponding $d$-modular overpartition $\mu$.
\begin{defi}For an overpartition $\lambda$ whose overlined parts are $\pm r,0 \pmod{d}$, nonoverlined parts are multiples of $d$, let $\mu$ denote the corresponding $d$-modular overpartition of $\lambda$. Then  the part  $\mu_i$ of $\mu$ is the number of elements in the $i$-th row of $d$-modular Ferrers diagram of $\lambda$, with the label according to the residue classes of $d$. The order of the parts in $d$-modular overpartitions is subject to that
\[1_{\overline{r}}<1_{\overline{d-r}}<1_{\overline{d}}<2_d<2_{\overline{r}}<2_{\overline{d-r}}<2_{\overline{d}}<\cdots.\]
The weight of the parts of $d$-modular overpartitions is given as follows:
\begin{itemize}
\item[(i)]the weight of part  $n_{\overline{r}}$ is $(n-1)d+r$;
\item[(ii)]the weight of part  $n_{\overline{d-r}}$ is $(n-1)d+d-r=nd-r$;
\item[(iii)]the weight of part  $n_{\overline{d}}$ is $(n-1)d+d=nd$;
\item[(iv)]the weight of part  $n_d$ is $(n-1)d+0=(n-1)d$.
\end{itemize}
\end{defi}
According to the overpartition in last example, $\mu=(4_{\overline{5}},4_7,4_7,3_{\overline{2}},1_{\overline{7}})$.

Now we define a matrix which describes the difference conditions between the parts of $d$-modular overpartition $\mu$ according to the label of the parts.

\begin{equation}
\overline{A'}_{d,r}=
\begin{blockarray}{ccccc}
& \overline{r} & \overline{d-r} & \overline{d} & d \\
\begin{block}{c(cccc)}
  \overline{r} & 1 & 1 & 2 & 0           \\
  \overline{d-r} & 1 & 1 & 2 & 0   \\
  \overline{d} & 1 & 1 & 2 & 0        \\
  d & 1 & 1 & 2 & 0            \\
\end{block}
\end{blockarray}
 \end{equation}

The rows and columns are indexed by $\overline{r}$, $\overline{d-r}$, $\overline{d}$, and $d$, so that, for example, $\overline{A'}_{d,r}(\overline{d}, \overline{d-r})=1$.

Now we focus on the overpartitions enumerated by $\overline{B}_{d,r}(n)$. For an overprtition $\lambda$ enumerated by $\overline{B}_{d,r}(n)$, we  consider the corresponding $d$-modular overpartition $\mu$, then we can see that the $d$-modular overpartition $\mu$ satisfies the following conditions (I):

\begin{itemize}
\item[(i)] the smallest part is odd and no $1_d$;
\item[(ii)] for $u, v \in \{\overline{r}, \overline{d-r}, \overline{d}, d\}$, if $\mu_{i+1}$ labeled  $u$  and $\mu_i$ labeled $v$ , then $\mu_{i+1}-\mu_i\geq \overline{A'}_{d,r}(u, v)$;
\item[(iii)] for $u, v \in \{\overline{r}, \overline{d-r}, \overline{d}, d\}$, if $\mu_{i+1}$ labeled  $u$  and $\mu_i$ labeled $v$ , then $\mu_{i+1}-\mu_i\equiv
\overline{A'}_{d,r}(u, v) \pmod{2}$.
\end{itemize}

That is to say if  $\mu$ is of the following form $\mu_1\geq \mu_2\geq \cdots\geq \mu_s$,
then (i)$\mu_i-\mu_{i+1}\geq 0$ and even, if $\mu_{i+1}$ labeled with $d$,
(ii)$\mu_i-\mu_{i+1}\geq 1$ and odd, if $\mu_{i+1}$ labeled with $\overline{r}$ or $\overline{d-r}$,
(iii)$\mu_i-\mu_{i+1}\geq 2$ and even, if $\mu_{i+1}$ labeled with $\overline{d}$.
(iv)no non-overlined $1$.

By employing the corresponding $d$-modular overpartition, the construction of the overpartitions enumerated by $\overline{B}_{d,r}(n)$ is much clearer. We also employ the $d$-modular partition to describe the overpartitions enumerated by $\overline{C}_{d,r}(n)$. For an overpartition $\nu$ enumerated by $\overline{C}_{d,r}(n)$ and the corresponding $d$-modular partition  $\xi$,  $\xi$ satisfies the following conditions:
\begin{itemize}
\item[(i)] the smallest part is even;
\item[(ii)] for $u, v \in \{r, d-r, \overline{d}, d\}$, if $\xi_{i+1}$ labeled  $u$  and $\xi_i$ labeled $v$ , then $\xi_{i+1}-\xi_i\geq \overline{A'}_{d,r}(u, v)$;
\item[(iii)] for $u, v \in \{r, d-r, \overline{d}, d\}$, if $\xi_{i+1}$ labeled  $u$  and $\xi_i$ labeled $v$ , then $\xi_{i+1}-\xi_i\equiv
\overline{A'}_{d,r}(u, v) \pmod{2}$.
\end{itemize}

Then we shall use the $d$-modular diagram and the corresponding $d$-modular partition to give the bijective proofs of Theorem \ref{thmB} and Theorem \ref{thmC}.

\section{The bijective proof of Theorem \ref{thmB}}

We shall give a bijective proof of Theorem \ref{thmB}, by using the $d$-modular Ferrers diagram and $d$-modular partition defined in the previous  section.

We just need to prove \eqref{BE}, which reads that
\[\overline{B}_{d,r}(n) =\overline{E}_{d,r}(n).\]

Recall that the overpartition  enumerated by $\overline{E}_{d,r}(n)$ is an overpartition with unrestricted parts congruence $0$ when modulo $2d$ and overlined parts congruent to $\pm r$ when modulo $d$. Then we can see these overpartitons also have $d$-modular Ferrers diagrams and corresponding $d$-modular overpartitions.

We let $\overline{E}_{d,r}$ denote the set of the overpartitions enumerated by $\overline{E}_{d,r}(n)$ for all $n\geq 0$.  Then any overpartitions $\pi$ in $\overline{E}_{d,r}$ corresponding to a $d$-modular partition is subjected to that,
\begin{itemize}
\item[1.]the non-overlined parts are all odd and greater than $2$.
 \item[2.]the overlined parts are all labeled with $d$ and $d-r$.
  \item[3.]the overlined parts with same label are distinct.
  \end{itemize}
For an overpartition  $\pi\in \overline{E}_{d,r}$, let $\chi$ denote the corresponding $d$-modula partition of $\pi$, then we decompose $\chi$ into  a partition triple $(\alpha,\beta,\gamma)$, where $\alpha$ is a partition with odd parts$\geq 3$ and all labeled with $d$, $\beta$ is a partition with distinct parts all labeled with $\overline{r}$, and $\gamma$ is a partition with distinct parts all labeled with $\overline{d-r}$.

We will construct a $d$-modula overpartition $\mu$ from the  partition triple
$(\alpha,\beta,\gamma)\in \overline{E}_{d,r}$, whose corresponding overpartition is $\lambda\in\overline{B}_{k,i}$. One will see that this procedure is invertible which makes the proof  bijective. By the definition of the weight of the parts in $d$-modular we can see that the construction also keeps the weight of the partition, that is, $|\pi|=|\chi|=|\alpha|+|\beta|+|\gamma|=|\mu|$.

We will insert the partitions $\beta$ and $\gamma$ into $\alpha$.
Let $l(*)$ denote the number of parts of a partition. We consider the largest part of $\beta$ and $\gamma$ and first insert the partition whose largest part is bigger. Then we discuss the following conditions.

\begin{itemize}
\item[Case 1.] $\beta_1\geq \gamma_1$.
We first insert the parts of $\beta$ into $\alpha$. We compare the largest part of $\beta$ and the number of parts of $\alpha$ and consider the following two cases.

\item[Case 1.1]$\beta_1> l(\alpha)$.

Step 1. In this step we shall get a $d$-modular partition $\alpha^{(1)}$ with exactly $\beta_1$ parts. The insertion will manipulate   the $d$-modular Ferrers diagram. We insert $\beta_1$ into $\alpha$ to get $\alpha^{(1)}$ by putting  $\beta_1$ as the first column of $\alpha^{(1)}$. Then  (i) for $1\leq i\leq l(\alpha)$, $\alpha^{(1)}_i=\alpha_i+1$, with the label   not changing; (ii) for $ l(\alpha)\leq i\leq \beta_1-1$, $\alpha^{(1)}_i=2$ labeled with $d$;(iii) $\alpha^{(1)}_{\beta_1}=1$ with label $\overline{r}$.  $|\alpha^{(1)}|=|\alpha|+|\beta_1|$.

Step 2. For $i\neq 1$, we insert any  other part $\beta_i$ of $\beta$ into $\alpha^{(i-1)}$ by adding a column to $\alpha^{(i-1)}$. More precisely, let $\alpha^{(i-1)}_1$,
$\alpha^{(i-1)}_2$,\ldots, $\alpha^{(i-1)}_{\beta_i-1}$ plus $1$ with the label   not changing, and the $\alpha^{(i-1)}_{\beta_i}$ just change the label from $d$ to  $\overline{r}$. Then after inserting all parts in $\beta$ we get a $d$-modular  partition $\alpha^{(l(\beta))}$. $|\alpha^{(l(\beta))}|=|\alpha|+|\beta|$.

Step 3. We insert the  parts of $\gamma$  into $\alpha^{(l(\beta))}$ successively. For a part $\gamma_i$, we  insert it into $\alpha^{(l(\beta)+i-1)}$ by adding one column. Let the parts $\alpha^{(l(\beta)+i)}_1=\alpha^{(l(\beta)+i-1)}_1+1$, $\alpha^{(l(\beta)+i)}_2=\alpha^{(l(\beta)+i-1)}_2+1$,\ldots, $\alpha^{(l(\beta)+i)}_{\gamma_i-1}=\alpha^{(l(\beta)+i-1)}_{\gamma_i-1}+1$ with labels   not changing. For $\alpha^{(l(\beta)+i)}_{\gamma_i}$, if $\alpha^{(l(\beta)+i-1)}$ is a part labeled with $d$, then we get  $\alpha^{(l(\beta)+i)}_{\gamma_i}=\alpha^{(l(\beta)+i-1)}_{\gamma_i}$ with label $\overline{d-r}$; if $\alpha^{(l(\beta)+i-1)}$ is a part labeled with $\overline{r}$, then we get  $\alpha^{(l(\beta)+i)}_{\gamma_i}=\alpha^{(l(\beta)+i-1)}_{\gamma_i}$ with label $\overline{d}$. Then we let $\mu=\alpha^{(l(\beta)+l(\gamma)}$ satisfing that  $|\mu|=|\lambda|+|\beta|+|\gamma|$.
\item[Case 1.2]$l(\alpha)\geq \beta_1$.

Step 1. From $i=1$ to $i=l(\beta)$, we insert $\beta_i$ successively. For any   part $\beta_i$ of $\beta$. We let $\alpha^{(i-1)}_1$,
$\alpha^{(i-1)}_2$,\ldots, $\alpha^{(i-1)}_{\beta_i-1}$ plus $1$ with the labels   not changing, and the $\alpha^{(i-1)}_{\beta_i}$ just change the label to be $\overline{r}$. Then after inserting all parts in $\beta$ we get the partition $\alpha^{(l(\beta))}$.

Step 2. From $i=1$ to $i=l(\gamma)$, we insert $\gamma_i$ into $\alpha^{(l(\beta))}$ successively. For a part $\gamma_i$, we  insert it into $\alpha^{(l(\beta)+i-1)}$ by letting the parts $\alpha^{(l(\beta)+i)}_1=\alpha^{(l(\beta)+i-1)}_1+1$, $\alpha^{(l(\beta)+i)}_2=\alpha^{(l(\beta)+i-1)}_2+1$,\ldots, $\alpha^{(l(\beta)+i)}_{\gamma_i-1}=\alpha^{(l(\beta)+i-1)}_{\gamma_i-1}+1$ with labels   not changing. For $\alpha^{(l(\beta)+i)}_{\gamma_i}$, if $\alpha^{(l(\beta)+i-1)}$ is a part labeled with $d$ then we get  $\alpha^{(l(\beta)+i)}_{\gamma_i}=\alpha^{(l(\beta)+i-1)}_{\gamma_i}$ with label $\overline{d-r}$; if $\alpha^{(l(\beta)+i-1)}$ is a part labeled with $\overline{r}$, then we get  $\alpha^{(l(\beta)+i)}_{\gamma_i}=\alpha^{(l(\beta)+i-1)}_{\gamma_i}$ with label $\overline{d}$. At last we get $\mu=\alpha^{l(\beta)+l(\gamma)}$ subjected to that $|\mu|=|\lambda|+|\beta|+|\gamma|$.

\item[Case 2.] $\beta_1<\gamma_1$, then we first insert $\gamma$ and then $\beta$. The steps are similar as Case 1, so we omit it here.
\end{itemize}
After the insertion we get a $d$-modula overpartition $\mu$ which satisfies the conditions (I) and the corresponding overpartition $\lambda \in\overline{B}_{d,r}$.

We give an example in Case 1.1. Let $d=7,r=2$, and $\pi=(42, 42,\overline{37},\overline{30},28, \overline{26},\overline{19},14,14,\overline{9},\overline{5},\overline{2})$. Then the corresponding $d$-modula partition triple is $((7_d,7_d,5_d,3_d,3_d),(6_{\overline{2}},5_{\overline{2}},2_{\overline{2}},1_{\overline{2}}),(4_{\overline{5}},3_{\overline{5}},1_{\overline{5}}))$.

We first display the $\alpha$ as follows:
\begin{center}
 \begin{picture}(100,100)
\put(0,15){$7$}\put(15,15){$7$}\put(30,15){$0$}

\put(0,30){$7$}\put(15,30){$7$}\put(30,30){$0$}

\put(15,45){$7$}
\put(0,45){$7$}
\put(30,45){$7$}\put(45,45){$7$}
\put(60,45){$0$}

\put(15,60){$7$}\put(0,60){$7$}
\put(30,60){$7$}
\put(45,60){$7$}\put(60,60){$7$}
\put(75,60){$7$}\put(90,60){$0$}

\put(0,75){$7$}
\put(15,75){$7$}
\put(30,75){$7$}
\put(45,75){$7$}
\put(60,75){$7$}
\put(75,75){$7$}
\put(90,75){$0$}
\put(0,-10){$d$-modular partition $\alpha$}
 \end{picture}
\end{center}
We can see the largest part of $\beta$ is $\beta_1=6_{\overline{2}}$,   greater than $\gamma_1=4_{\overline{5}}$. So we consider the Case 1. Since $l(\alpha)=5<\beta_1$, so we consider the case 1.1,namely, we should insert $6_{\overline{2}}$ into $\alpha$ to get
\begin{center}
 \begin{picture}(150,100)
\put(30,30){$\overline{2}$}

\put(30,45){$7$}\put(45,45){$7$}\put(60,45){$7$}\put(75,45){$0$}

\put(30,60){$7$}\put(45,60){$7$}\put(60,60){$7$}\put(75,60){$0$}

\put(45,75){$7$}\put(30,75){$7$}\put(60,75){$7$}\put(75,75){$7$}
\put(90,75){$7$}\put(105,75){$0$}

\put(45,90){$7$}\put(30,90){$7$}
\put(60,90){$7$}
\put(75,90){$7$}\put(90,90){$7$}
\put(105,90){$7$}\put(120,90){$7$}\put(135,90){$0$}

\put(30,105){$7$}
\put(45,105){$7$}
\put(60,105){$7$}
\put(75,105){$7$}
\put(90,105){$7$}
\put(105,105){$7$}
\put(120,105){$7$}
\put(135,105){$0$}
\put(20,10){$\alpha^{(1)}$: by inserting $6_{\overline{2}}$ into $\alpha$}

 \end{picture}
\end{center}
Now we successively insert other parts of $\beta$ into $\alpha^{(1)}$ to get $\alpha^{(4)}$ as follows.
\begin{center}
 \begin{picture}(200,100)
\put(40,0){$\overline{2}$}

\put(40,15){$7$}\put(55,15){$7$}\put(70,15){$7$}\put(85,15){$\overline{2}$}

\put(40,30){$7$}\put(55,30){$7$}\put(70,30){$7$}\put(85,30){$7$}\put(100,30){$0$}

\put(55,45){$7$}
\put(40,45){$7$}
\put(70,45){$7$}\put(85,45){$7$}
\put(100,45){$7$}
\put(115,45){$7$}
\put(130,45){$0$}

\put(55,60){$7$}\put(40,60){$7$}
\put(70,60){$7$}
\put(85,60){$7$}\put(100,60){$7$}
\put(115,60){$7$}\put(130,60){$7$}\put(145,60){$7$}\put(160,60){$\overline{2}$}

\put(40,75){$7$}
\put(55,75){$7$}
\put(70,75){$7$}
\put(85,75){$7$}
\put(100,75){$7$}
\put(115,75){$7$}
\put(130,75){$7$}
\put(145,75){$7$}
\put(160,75){$7$}
\put(175,75){$\overline{2}$}
 \end{picture}
\end{center}
At last, we insert the parts of $\gamma$ into the above $d$-modular diagram of $\alpha^{(4)}$ to get $\alpha^{(7)}$ as follows.
\begin{center}
 \begin{picture}(250,100)
\put(40,0){$\overline{2}$}

\put(40,15){$7$}\put(55,15){$7$}\put(70,15){$7$}\put(85,15){$\overline{2}$}

\put(40,30){$7$}\put(55,30){$7$}\put(70,30){$7$}\put(85,30){$7$}\put(100,30){$\overline{5}$}

\put(55,45){$7$}
\put(40,45){$7$}
\put(70,45){$7$}\put(85,45){$7$}
\put(100,45){$7$}
\put(115,45){$7$}
\put(130,45){$7$}
\put(145,45){$\overline{5}$}

\put(55,60){$7$}\put(40,60){$7$}
\put(70,60){$7$}
\put(85,60){$7$}\put(100,60){$7$}
\put(115,60){$7$}\put(130,60){$7$}\put(145,60){$7$}
\put(160,60){$7$}\put(175,60){$7$}\put(190,60){$\overline{2}$}

\put(40,75){$7$}
\put(55,75){$7$}
\put(70,75){$7$}
\put(85,75){$7$}
\put(100,75){$7$}
\put(115,75){$7$}
\put(130,75){$7$}
\put(145,75){$7$}
\put(160,75){$7$}
\put(175,75){$7$}
\put(190,75){$7$}
\put(205,75){$\overline{7}$}

 \end{picture}
\end{center}

The inverse map is a subtraction, one can get it easily, so we omit it here.

\section{Construction of the generating function of $\overline{C}_{d,r}(n)$}

We shall give an alternative form of the generating function of $\overline{C}_{d,r}(n)$ in the combinatorial way. Then by a $_3\phi_2$ transformation we can prove that this generating function form  equals  the generating function of $\overline{C}_{d,r}(n)$ in Theorem \ref{thmC}.

For an overpartition $\lambda$ enumerated by $\overline{C}_{d,r}(n)$, we also use the $d$-modular Ferrers diagram. If the smallest part of $\lambda$ is at least $d$, then to construct any  overpartition $\lambda$ with exactly $m$ parts, we start from the partition with $m$ non-overlined $d$. Display by the $d$-modular Ferrers diagram, the overpartition has two columns with length $m$, where the left one has $m$ $d$'s and the right one has $m$ $0$'s. Then we construct all overpartitions enumerated by $\overline{C}_{d,r}(n)$ with exactly $m$ parts by inserting a partition triple $(\alpha,\beta,\gamma)$, where $\alpha$ is a partition with parts $\leq 2md$ and all parts are even multiple of $d$, $\beta$ is an overpartition  with parts $\leq (m-1)d+r$ are all congruent to  $r$ when modulo $d$, and $\gamma$ is an overpartition with  all parts $\leq md-r$ and congruent to $d-r$ when modulo $d$.

\noindent The construction of the overpartitions enumerated by $\overline{C}_{d,r}(n)$.
\begin{itemize}
\item[Step 1]Put a column with  $m$ $d$'s and a column with $m$ $0$'s. The weight is $md$. For example, we construct an overpartition $\lambda\in \overline{C}_{7,2}$ with $6$ pats, whose weight is $42$.
\begin{center}
 \begin{picture}(100,90)

\put(40,0){$7$}\put(40,15){$7$}\put(40,30){$7$}\put(40,45){$7$}
\put(40,60){$7$}\put(40,75){$7$}
\put(55,0){$0$}\put(55,15){$0$}\put(55,30){$0$}\put(55,45){$0$}
\put(55,60){$0$}\put(55,75){$0$}
\end{picture}
\end{center}

\item[Step 2.] Recall that $\alpha$ is a partition with all parts are multiples of $2d$ and the largest parts $\leq 2md$. By inserting all parts in  $\alpha$ into $\lambda^{(0)}=(d_d,d_d,\ldots,d_d)$ successively to get $\lambda^{(1)}$. For a part $\alpha_i=2sd$ we insert it by enlarging the first $s$ parts $2d$.
 Then after inserting all parts in $\alpha$ we get a $d$-modular Ferrers diagram with $m$ rows subject to that each row with odd number of $d$'s and one $0$.
 The weight of $\lambda^{(1)}$ is $|\lambda^{(1)}|=md+|\alpha|$.

 We go on with the example given in the Step 1. We will insert the partition $(70,42,42,14)$ into the $d$-modular Ferrers diagram.
 \begin{center}
 \begin{picture}(150,90)
  \put(10,0){$7$}\put(10,15){$7$}\put(10,30){$7$}\put(10,45){$7$}
\put(10,60){$7$}\put(10,75){$7$}

  \put(25,0){$0$}\put(25,15){$7$}\put(25,30){$7$}\put(25,45){$7$}
\put(25,60){$7$}\put(25,75){$7$}

\put(40,15){$7$}\put(40,30){$7$}\put(40,45){$7$}
\put(40,60){$7$}\put(40,75){$7$}
\put(55,15){$0$}\put(55,30){$0$}\put(55,45){$7$}
\put(55,60){$7$}\put(55,75){$7$}

\put(70,45){$7$}
\put(70,60){$7$}\put(70,75){$7$}

\put(85,45){$0$}
\put(85,60){$0$}\put(85,75){$7$}

\put(100,75){$7$} \put(115,75){$0$}

\end{picture}
\end{center}
\item[Step 3.]In this step we insert the parts of $\beta$ into $\lambda^{(1)}$ which are congruent to $r$ when modulo $d$. For each part $\overline{(k-1)d+r}$ in $\beta$, we insert it into $\lambda^{(1)}$ by letting the first $k-1$ parts of enlarge $d$ and the $k$th part enlarge $r$ and be an overlined part to get $\lambda^{(3)}$. The weight of $\lambda^{(3)}$ is $|\lambda^{(3)}|=|\lambda^{(2)}|+|\beta|$.

    Then suppose $\beta=(23,9,2)$, after inserting $\beta$ into $\lambda^{(1)}$, we get the following $d$-modular Ferrers diagram.
 \begin{center}
 \begin{picture}(150,90)
  \put(10,0){$7$}\put(10,15){$7$}\put(10,30){$7$}\put(10,45){$7$}
\put(10,60){$7$}\put(10,75){$7$}

  \put(25,0){$0$}\put(25,15){$7$}\put(25,30){$7$}\put(25,45){$7$}
\put(25,60){$7$}\put(25,75){$7$}

\put(40,15){$7$}\put(40,30){$7$}\put(40,45){$7$}
\put(40,60){$7$}\put(40,75){$7$}
\put(55,15){$0$}\put(55,30){$\overline{2}$}\put(55,45){$7$}
\put(55,60){$7$}\put(55,75){$7$}

\put(70,45){$7$}
\put(70,60){$7$}\put(70,75){$7$}

\put(85,45){$7$}
\put(85,60){$7$}\put(85,75){$7$}
\put(100,60){$\overline{2}$}

\put(100,75){$7$} \put(115,75){$7$}
\put(130,75){$7$}\put(145,75){$\overline{2}$}
\end{picture}
\end{center}
\item[Step 4.]In this step we insert the parts of $\gamma$ into $\lambda^{(2)}$ which are congruent to $d-r$ when modulo $d$. For each part $\overline{kd-r}$ in $\gamma$, insert it into $\lambda^{(2)}$ by letting the first $k-1$ parts enlarge $d$ and the $k$th part enlarge $d-r$ with overlined. We can see that if there is no $\overline{(k-1)d+r}$ in $\beta$, after this insertion, the $k$-th part becomes an overlined part congruent to $d-r$ modulo $d$; otherwise the $k$-th part become an overlined part congruent to $0$  modulo $d$.

    Suppose $\gamma=(26,19,5)$, we insert $\gamma$ into $\lambda^{(2)}$ to get the following
 \begin{center}
 \begin{picture}(150,90)
  \put(10,0){$7$}\put(10,15){$7$}\put(10,30){$7$}\put(10,45){$7$}
\put(10,60){$7$}\put(10,75){$7$}

  \put(25,0){$0$}\put(25,15){$7$}\put(25,30){$7$}\put(25,45){$7$}
\put(25,60){$7$}\put(25,75){$7$}

\put(40,15){$7$}\put(40,30){$7$}\put(40,45){$7$}
\put(40,60){$7$}\put(40,75){$7$}
\put(55,15){$0$}\put(55,30){$\overline{7}$}\put(55,45){$7$}
\put(55,60){$7$}\put(55,75){$7$}

\put(70,45){$7$}\put(85,45){$7$}\put(100,45){$\overline{5}$}
\put(70,60){$7$}\put(70,75){$7$}

\put(85,45){$7$}
\put(85,60){$7$}\put(85,75){$7$}
\put(100,60){$7$}
\put(115,60){$7$}
\put(130,60){$\overline{2}$}

\put(100,75){$7$} \put(115,75){$7$}
\put(130,75){$7$}\put(145,75){$7$}
\put(160,75){$7$}\put(175,75){$\overline{7}$}
\end{picture}
\end{center}
\end{itemize}
By this insertion we can get all the overpartition with $m$ parts enumerated by $\overline{C}_{d,r}(n)$.

By this construction we can get the generating function of $\overline{C}_{d,r}(m,n)$, which is the number of overpartitions enumerated by $\overline{C}_{d,r}(n)$ with exactly $m$ parts.
\begin{thm}
The generating function of $\overline{C}_{d,r}(m,n)$ is as follows.

\begin{equation}\label{cmn}\sum_{n\geq 0}\overline{C}_{d,r}(m,n)x^mq^n=\sum_{m\geq 0}\frac{x^mq^{dm}(-q^r,-q^{d-r};q^d)_m}{(q^{2d};q^{2d})_m}.\end{equation}

\end{thm}
We shall prove that by letting $x=1$. Then, the right hand side of \eqref{cmn}  equals the right hand side of \eqref{bc}, that is
\begin{equation}\label{cc}\sum_{n\geq 0}\frac{q^{dn}(-q^r,-q^{d-r};q^d)_n}{(q^{2d};q^{2d})_n}=\frac{(-q^r,-q^{d-r};q^d)_\infty}{(q^{2d};q^{2d})_\infty}\sum_{n\geq 0}\frac{(-q^d;q^d)_nq^{\frac{dn(n+1)}{2}}}{(-q^r,-q^{d-r};q^d)_{n+1}}.
\end{equation}

We employ the  following $_3\phi_2$ transformation (which is found in an equivalent form as equation (III.10) in [9])
\begin{equation}\sum_{n\geq 0}\frac{(\frac{aq}{bc},d,e;q)_n}{(q,\frac{aq}{b},\frac{aq}{c};q)_n}(\frac{aq}{de})^n
=\frac{(\frac{aq}{d},\frac{aq}{e},\frac{aq}{bc};q)_\infty}{(\frac{aq}{b},\frac{aq}{c},\frac{aq}{de};q)_\infty}
\sum_{n\geq0}\frac{(\frac{aq}{de},b,c;q)_n}{(q,\frac{aq}{d},\frac{aq}{e};q)_n}(\frac{aq}{bc})^n.
\end{equation}
 Let $b\rightarrow\infty$ and set $q=q^d,\ a=q^d,\ c=-q^d,\ d=-q^r,\ e=-q^{d-r}$, we get the identity \eqref{cc}.

\vspace{0.5cm}
 \noindent{\bf Acknowledgments.}
This work was supported by the National Science Foundation of China (Nos.1140149, 11501089, 11501408).

\end{document}